%% file: Elasto_3D.tex
\documentclass{article}
\usepackage{stmaryrd}
\usepackage{mathrsfs}
\usepackage{amsmath}
\usepackage{amssymb}
\usepackage{amsfonts}
\usepackage{pstricks}
\usepackage{graphicx, color, epsfig, enumerate}

\newcommand{\om}{\omega}
\newcommand{\ga}{\gamma}

\newcommand{\la}{\lambda}
\newcommand{\q}{\tilde}
\newcommand{\al}{\alpha}

\newcommand{\f}{\frac}

\newtheorem{theorem}{Theorem}[section]
\newtheorem{proposition}[theorem]{Proposition}

\numberwithin{equation}{section}

\begin{document}
\title{{A Novel Approach to Elastodynamics:}\\
{II. The Three-Dimensional Case}}
\author{
{A. S. Fokas${}^a$\footnote{T.Fokas@damtp.cam.ac.uk},\  \
D. Yang${}^{a,}{}^{b}$\footnote{yangd04@mails.tsinghua.edu.cn}}\\
\\
{\small ${}^a$ Department of Applied Mathematics and Theoretical Physics,}\\
{\small University of Cambridge, Cambridge CB3 0WA, UK}\\
%{\small ~~~~~~~~~~D.Yang@damtp.cam.ac.uk}\\
\\
{\small ${}^b$ Department of Mathematical Sciences, Tsinghua
University,}\\
{\small Beijing 100084, P. R. China}\\
}
\maketitle

\begin{abstract}
A new approach was recently introduced by the authors for constructing analytic solutions of
the linear PDEs describing elastodynamics. Here, this approach is applied to
the case of a homogeneous isotropic half-space body
satisfying arbitrary initial conditions and Lamb's boundary
conditions. A particular case of this problem, namely the case of
homogeneous initial conditions and normal point load boundary conditions,
was first solved by Lamb using the Fourier-Laplace transform.
The general problem solved here can also be analysed via the Fourier transform, but in this case, the solution
representation involves transforms of \textit{unknown} boundary
values; this necessitates the formulation and solution of a
cumbersome auxiliary problem, which expresses the unknown boundary
values in terms of the Laplace transform of the given boundary
data. The new approach, which is applicable to arbitrary
initial and boundary conditions, bypasses the above auxiliary
problem and expresses the solutions directly in terms of the given
initial and boundary conditions.
%The initial boundary value problem(IBVP) for 2D-elastodynamics in
%the half space is studied through a unified approach introduced in
%1990s by A. S. Fokas\cite{ASF}. This IBVP is a long-standing problem
%since 1900s. And the unified approach for solving this IBVP is based
%on an intensive complex analysis of the so-called global relations.
%Compared with the known Fourier-Laplace transform method first
%employed to this IBVP by Lamb, our method is applied not only to the
%general boundary conditions but also to the general initial
%conditions; and the final exact solution is expressed in an elegant
%and closed form.
\end{abstract}

\section{Introduction}
The problem considered in this paper has a long and illustrious
history, which begins with the classic works of Sir Horace Lamb in
1904 \cite{HL}. In \cite{HL}, Lamb treated four basic problems, the
so called Lamb's problems, which are formulated in either two
or three dimensions. %, the surface normal line and
%point load, and the buried line and point load.

These problems have been studied by several authors, see for example
\cite{JM, KLJ, point_load0, point_load05, point_load1, point_load2,
point_load4, point_load5, line_load1, line_load16, line_load2,
Payton}.

A new approach to elastodynamics was introduced in \cite{ASF_Yang}.
This approach is based on the unified method for solving linear and
integrable nonlinear PDEs introduced by one of the authors in
\cite{ASF1997}. In \cite{ASF_Yang}, the new approach was applied to
the Lamb's problem in two dimensions. Here, we study
three-dimensional problems. In particular, we consider arbitrary
initial conditions and general stress boundary conditions, including
normal point load, tangential point load, and mixed point load. We
refer to the latter three stress conditions as Lamb's boundary
conditions.

Most studies of Lamb's problems were based on the Helmholtz
decomposition and on the use of the Laplace transform in time.
Helmholtz decomposition has the advantage of decomposing the P-waves
and S-waves. However, it has the disadvantage of introducing higher
order derivatives to the boundary conditions. The use of Laplace
transform, essentially restricts the problem to the case of
homogeneous initial conditions.

A crucial role of the method of Fokas \cite{ASF1997}\cite{ASF} is
played by the so-called \textit{global relations}, which are
algebraic relations coupling appropriate transforms of the unknown
boundary values with transforms of the given data. For linear PDEs
the method of Fokas uses three novel steps \cite{ASF}-\cite{SSF}: 1.
Derive a representation for the solution in terms of an integral
involving a contour in the complex Fourier plane. This
representation is not yet effective, because in addition to
transforms of the given initial and boundary data, it also contains
transforms of unknown boundary values. 2. Analyse certain
transformations in the complex Fourier plane which leave invariant
the transforms of the unknown boundary values. 3. Eliminate the
transforms of the unknown boundary values, by combining the results
of the first two steps and by employing Cauchy's theorem (or more
precisely Jordan's lemma).

This paper is organised as follows. In section 2, we recall the
governing equations of elastodynamics and derive the global
relations. In section 3 we implement step 1. In section 4 we
implement steps 2 and 3. In section 5, by employing the general
representations derived in section 4, we analyse Lamb's problems.
These results are further discussed in section 6.

For certain complicated boundary value problems it seems that it is
not possible to eliminate from the integral representation of the
solution the transforms of the unknown boundary values. However, for
some of these problems, by using the global relations, one can
derive expressions for the Laplace transforms of the
\textit{unknown} boundary values in terms of the given data
\cite{Ashton_ASF}. A summary of this \textit{less} effective
approach is presented in the Appendix. It is interesting that for
the \textit{particular case of zero initial conditions}, the
formulae presented in the Appendix reduce to the formulae first
derived in the classic works of Lamb.
%We also give the exact solution for initial value problems of planar
%elastodynamics in Appendix\ref{IVP}.
\section{Governing Equations and Global Relations}
The transient problem for three dimensional elastodynamics in the half space with
Lamb's boundary conditions is defined as follows: Let
$u=u(x,y,z,t),$ $v=v(x,y,z,t),$ $w=w(x,y,z,t),$ denote the displacements of a
homogeneous isotropic half space body. The governing equations of
motion without external body forces, are the Lam\'e-Navier
equations:
\begin{subequations}\label{3D}
\begin{align}
&(\la+2\mu)u_{xx}+(\la+\mu)(v_{xy}+w_{xz})+\mu (u_{yy}+u_{zz})-\rho u_{tt}=0,\\
&(\la+2\mu)v_{yy}+(\la+\mu)(w_{yz}+u_{yx})+\mu (v_{zz}+v_{xx})-\rho v_{tt}=0,\\
&(\la+2\mu)w_{zz}+(\la+\mu)(u_{zx}+v_{zy})+\mu (w_{xx}+w_{yy})-\rho w_{tt}=0,
\end{align}
\begin{flushright}
$~-\infty<x<\infty,~-\infty<y<\infty, ~~z>0,~~ t>0,$
\end{flushright}
\end{subequations}
where $\la,\mu$ are the Lam\'e constants and $\rho$ denotes the density of
the material, which without loss of generality is normalized to unity, i.e., $\rho\equiv1.$
Let the initial conditions be denoted by
\begin{subequations}\label{ic}
\begin{align}
&u(x,y,z,0)=u_0(x,y,z),~u_t(x,y,z,0)=u_1(x,y,z),\\
&v(x,y,z,0)=v_0(x,y,z),~v_t(x,y,z,0)=v_1(x,y,z),\\
&w(x,y,z,0)=w_0(x,y,z),~w_t(x,y,z,0)=w_1(x,y,z),
\end{align}
\begin{flushright}$~-\infty<x<\infty,~-\infty<y<\infty, ~~z>0. $\end{flushright}
\end{subequations}

Let the stress boundary conditions be denoted by
\begin{subequations}\label{3D_l_bv}
 \begin{align}
(u_z+w_x)(x,y,0,t)=g_1(x,y,t),\\
(v_z+w_y)(x,y,0,t)=g_2(x,y,t),\\
\Big(w_z+\f{\la}{\la+2\mu}(u_x+v_y)\Big)(x,y,0,t)=g_3(x,y,t),
\end{align}
\begin{flushright}$~-\infty<x<\infty,~-\infty<y<\infty, ~~t>0. $\end{flushright}
\end{subequations}
For a tangential point load, the functions $g_1$, $g_2$ and $g_3$ are given by
\begin{equation}
g_1(x,y,t)=\sigma_0\delta(x,y)X(t)/\mu,~g_2(x,y,t)=0,~g_3(x,y,t)=0;
\end{equation}
for a normal point load,
\begin{equation}
g_1(x,y,t)=0,~g_2(x,t)=0,~g_3(x,y,t)=\sigma_0\delta(x,y)Y(t)/(\la+\mu);
\end{equation}
for a moving normal line load with a constant speed $C$ along the $x-axis$,
\begin{equation}
g_1(x,y,t)=0,~g_2(x,y,t)=0,~g_3(x,y,t)=\sigma_0\delta(x-Ct,y)/(\la+\mu).
\end{equation}
Here, $\sigma_0$ is a constant, $\delta(x,y)$ denotes the Dirac-$\delta$ function, $X(t)$ and $Y(t)$ are functions which depend only on $t$.
\paragraph{Notations}
Hat, ``$\wedge$'', will denote the three dimensional Fourier transform with respect to
$x$, $y$ and $z$, whereas tilde, ``$\sim$'', will denote the two dimensional Fourier
transform with respect to $x$ and $y$. In particular,
\begin{subequations}\label{Fourier}
\begin{align}
&\hat{u}(k,l,m,t):=\int_{-\infty}^{\infty}dx\int_{-\infty}^{\infty}dy\int_0^{\infty}dz~
e^{-ikx-ily-imz}u(x,y,z,t),\\
&\hat{v}(k,l,m,t):=\int_{-\infty}^{\infty}dx\int_{-\infty}^{\infty}dy\int_0^{\infty}dz~
e^{-ikx-ily-imz}v(x,y,z,t),\\
&\hat{w}(k,l,m,t):=\int_{-\infty}^{\infty}dx\int_{-\infty}^{\infty}dy\int_0^{\infty}dz~
e^{-ikx-ily-imz}w(x,y,z,t),\\
&\hat{u}_j(k,l,m,t):=\int_{-\infty}^{\infty}dx\int_{-\infty}^{\infty}dy\int_0^{\infty}dz~
e^{-ikx-ily-imz}u_j(x,y,z,t),\\
&\hat{v}_j(k,l,m,t):=\int_{-\infty}^{\infty}dx\int_{-\infty}^{\infty}dy\int_0^{\infty}dz~
e^{-ikx-ily-imz}v_j(x,y,z,t),\\
&\hat{w}_j(k,l,m,t):=\int_{-\infty}^{\infty}dx\int_{-\infty}^{\infty}dy\int_0^{\infty}dz~
e^{-ikx-ily-imz}w_j(x,y,z,t),
\end{align}
\begin{flushright}$k,l\in\mathbb{R},~m\in\mathbb{C}^-,~t>0, j=0,1,$\end{flushright}
\end{subequations}
where $\mathbb{C}^-$ denotes the lower half complex $m$-plane.

Furthermore,
\begin{subequations}\label{xF}
\begin{align}
\tilde{u}(k,l,t)&:=\int_{-\infty}^{\infty}dx \int_{-\infty}^{\infty}dy~~
e^{-ikx-ily}u(x,y,0,t),\\
\tilde{v}(k,l,t)&:=\int_{-\infty}^{\infty}dx \int_{-\infty}^{\infty}dy~~
e^{-ikx-ily}v(x,y,0,t),\\
\tilde{w}(k,l,t)&:=\int_{-\infty}^{\infty}dx \int_{-\infty}^{\infty}dy~~
e^{-ikx-ily}w(x,y,0,t),  ~~ k,l\in\mathbb{R},t>0.
\end{align}
\end{subequations}
\begin{equation}\label{gtilde}
\q{g}_j(k,l,t):=\int_{-\infty}^{\infty}dx \int_{-\infty}^{\infty}dy~~e^{-ikx-ily} g_j(x,y,t),~~k,l\in\mathbb{R},~t>0,~j=1,2,3.
\end{equation}

We emphasise that since $x,y\in \mathbb{R}$, the $x,y-$Fourier transform
is well-defined only for $k,l\in\mathbb{R};$ on the other hand, since
$0<z<\infty$,  the $z-$Fourier transform is well-defined for $m$ %at least
in the lower half complex $m-$plane.

By applying the three dimensional Fourier transform to equations \eqref{3D} and by using in the resulting equations the boundary conditions \eqref{3D_l_bv}, we obtain the following equations:
\begin{subequations}\label{pre_g}
\begin{align}
\begin{split}
-\big[(\la+2\mu)k^2+\mu l^2+\mu m^2\big]\hat{u}-&(\la+\mu)kl\hat{v}-(\la+\mu)km\hat{w}\\
&-\mu\q{g}_1-\mu im \q{u}-\la ik\q{w}=\hat{u}_{tt},\\
\end{split}\\
\begin{split}
-\big[(\la+2\mu)l^2+\mu k^2+\mu m^2\big]\hat{v}-&(\la+\mu)kl\hat{u}-(\la+\mu)lm\hat{w}\\
&-\mu\q{g}_2-\mu im \q{v}-\la il\q{w}=\hat{v}_{tt},\\
\end{split}\\
\begin{split}
-\big[(\la+2\mu)m^2&+\mu k^2+\mu l^2\big]\hat{w}-(\la+\mu)km\hat{u}-(\la+\mu)lm\hat{v}\\
&-(\la+2\mu)\q{g}_3-(\la+2\mu)im \q{w}-\mu ik\q{u}-\mu
il\q{v}=\hat{w}_{tt},\\
\end{split}
\end{align}
\begin{flushright}$k,l\in\mathbb{R},~m\in\mathbb{C}^-.$\end{flushright}
\end{subequations}
Introducing the notations
\begin{subequations}
\begin{align}
 &P(k,l,m,t)=k\hat{u}(k,l,m,t)+l\hat{v}(k,l,m,t)+m\hat{w}(k,l,m,t),\\
 &Q(k,l,m,t)=m\hat{v}(k,l,m,t)-l\hat{w}(k,l,m,t),\\
 &R(k,l,m,t)=k\hat{w}(k,l,m,t)-m\hat{u}(k,l,m,t),%\\
 %&Q_3(k,l,m,t)=l\hat{u}(k,l,m,t)-k\hat{v}(k,l,m,t),
\end{align}
\end{subequations}
equations \eqref{pre_g} become
\begin{subequations}\label{ODEs}
\begin{align}
\label{ODE1} &P_{tt}+(\la+2\mu)(k^2+l^2+m^2)P=F_P,\\
\label{ODE2} &Q_{tt}+\mu(k^2+l^2+m^2)Q=F_{Q},\\
\label{ODE3} &R_{tt}+\mu(k^2+l^2+m^2)R=F_{R},~~k,l\in\mathbb{R},~m\in\mathbb{C}^-,
%\\  \label{ODE4} &Q_{3tt}+\mu(k^2+l^2+m^2)Q_3=F_{Q3},
\end{align}
\end{subequations}
where the functions $F_P(k,l,m,t)$, $F_Q(k,l,m,t)$ and $F_R(k,l,m,t)$ are defined as follows:
\begin{subequations}\label{Bpq}
\begin{align}
\label{Bp}
\begin{split}
&F_P=-(\mu k\q{g}_1+\mu l \q{g}_2+i\la k^2\q{w}+i\la l^2\q{w})-m (2i\mu k \q{u}+2i\mu l\q{v}\\
&~~~~~~~+(\la+2\mu)\q{g}_3)-i m^2(\la+2\mu)\q{v},\\
\end{split}\\
\label{Bq1}
&F_{Q}=((\la+2\mu)l\q{g}_3+i\mu kl \q{u}+i\mu l^2 \q{v})+m(-\mu\q{g}_2+2i\mu l\q{w})-m^2 i\mu \q{v},\\
\label{Bq2}
&F_{R}=-((\la+2\mu)k\q{g}_3+i\mu k^2 \q{u}+i\mu kl
\q{v})+m(\mu\q{g}_1-2i\mu k\q{w})+m^2 i\mu \q{u}.
\end{align}
\end{subequations}
Solving equations \eqref{ODEs} for $\{P,Q,R\}$, we obtain the following expressions:
\begin{subequations}\label{global_r}
\begin{align}
\label{global_r1}
\begin{split}
P=&\f{1}{-2i\om_0}\Big(e^{-i\om_0t}\int_0^t
e^{i\om_0s}F_P(k,l,m,s)ds-e^{i\om_0t}\int_0^t e^{-i\om_0s}F_P(k,l,m,s)ds\Big) \\
&+\Big(\f{1}{2}P_0+\f{i}{2}\f{P_1}{\om_0}\Big)e^{-i\om_0t}+\Big(\f{1}{2}P_0-\f{i}{2}\f{P_1}{\om_0}\Big)e^{i\om_0t},
\end{split}\\
\label{global_r2}
\begin{split}
Q=&\f{1}{-2i\om_1}\Big(e^{-i\om_1t}\int_0^t
e^{i\om_1s}F_{Q}(k,l,m,s)ds-e^{i\om_1t}\int_0^t e^{-i\om_1s}F_{Q}(k,l,m,s)ds\Big)\\
&+\Big(\f{1}{2}Q_0+\f{i}{2}\f{Q_1}{\om_1}\Big)e^{-i\om_1t}+\Big(\f{1}{2}Q_0-\f{i}{2}\f{Q_1}{\om_1}\Big)e^{i\om_1t},
\end{split}\\
\label{global_r3}
\begin{split}
R=&\f{1}{-2i\om_1}\Big(e^{-i\om_1t}\int_0^t
e^{i\om_1s}F_{R}(k,l,m,s)ds-e^{i\om_1t}\int_0^t e^{-i\om_1s}F_{R}(k,l,m,s)ds\Big)\\
&+\Big(\f{1}{2}R_0+\f{i}{2}\f{R_1}{\om_1}\Big)e^{-i\om_1t}+\Big(\f{1}{2}R_0-\f{i}{2}\f{R_1}{\om_1}\Big)e^{i\om_1t},
~~~~k,l\in\mathbb{R},~m\in\mathbb{C}^-,
\end{split}
\end{align}
\end{subequations}
where the dispersion relations $\om_1$ and $\om_2$ are given by
\begin{equation}\label{dispersion}
\om_1^2=(\la+2\mu)(k^2+l^2+m^2),~~\om_2^2=\mu(k^2+l^2+m^2),~~k,l\in\mathbb{R},~m\in\mathbb{C},
\end{equation}
and the \textit{known} functions $P_j(k,l,m),~Q_j(k,l,m),~R_j(k,l,m)$, $j=0,1$ are given in terms of the
initial conditions by
\begin{subequations}\label{initial}
 \begin{align}
&P_0=k\hat{u}_0+l\hat{v}_0+m\hat{w}_0,~~P_1=k\hat{u}_1+l\hat{v}_1+m\hat{w}_1,\\
&Q_0=m\hat{v}_0-l\hat{w}_0,~Q_1=m\hat{v}_1-l\hat{w}_1,\\
&R_0=k\hat{w}_0-m\hat{v}_0,~R_1=k\hat{w}_1-m\hat{v}_1,
 \end{align}
\end{subequations}
with $\hat{u}_j,~\hat{v}_j,~\hat{w}_j,~j=0,1,$ defined in equations \eqref{Fourier}.

In the following, we take
\begin{equation}\label{dispersion_choice}
\om_1=\sqrt{\la+2\mu}(k^2+l^2+m^2)^{\f{1}{2}},~~\om_2=\sqrt{\mu}(k^2+l^2+m^2)^{\f{1}{2}}.
\end{equation}
The function $(k^2+l^2+m^2)^{\f{1}{2}}$ has the branch points $\pm
i\sqrt{k^2+l^2}$ in the complex $m-$plane; we connect these two branch points by
a branch cut and we fix a branch in the cut plane by the requirement
that,
$$(k^2+l^2+m^2)^{\f{1}{2}}\sim~ m+O\Big(\f{1}{m}\Big),~~~as~~m\rightarrow\infty.$$
%For $k=l=0$, a direct result of eqs.(\ref{pre_g1})\&(\ref{pre_g2})
%when considering the initial conditions\eqref{ic1}\&\eqref{ic2}
%gives,
%\begin{subequations}\label{recall}
%\begin{align}
%\hat{u}(0,0,t)=\hat{u}_0(0,0)+t\hat{u}_1(0,0)-\mu
%\int_0^td\tau\int_0^{\tau}\q{g}_1(0,s)ds,\\
%\hat{v}(0,0,t)=\hat{v}_0(0,0)+t\hat{v}_1(0,0)-(\la+2\mu)
%\int_0^td\tau\int_0^{\tau}\q{g}_2(0,s)ds.
%\end{align}
%\end{subequations}

Let $u^{(j)\pm},~v^{(j)\pm},~w^{(j)\pm},~U^{(j)},~V^{(j)},~W^{(j)},~j=1,2,$ denote the
following \textit{unknown} functions:
\begin{subequations}\label{BigU}
 \begin{align}
&u^{(j)\pm}(k,l,m,t)=\int_0^t e^{\pm i\om_j s}\q{u}(k,l,s)ds,\\
&v^{(j)\pm}(k,l,m,t)=\int_0^t e^{\pm i\om_j s}\q{v}(k,l,s)ds,\\
&w^{(j)\pm}(k,l,m,t)=\int_0^t e^{\pm i\om_j s}\q{w}(k,l,s)ds,\\
\label{defU1}
&U^{(j)}(k,l,m,t)=\f{1}{2\om_j}\big(e^{-i\om_j t}u^{(j)+}(k,l,m,t)-e^{i\om_j t}u^{(j)-}(k,l,m,t)\big),\\
&V^{(j)}(k,l,m,t)=\f{1}{2\om_j}\big(e^{-i\om_j
t}v^{(j)+}(k,l,m,t)-e^{i\om_j t}v^{(j)-}(k,l,m,t)\big),\\
\label{defW1}
&W^{(j)}(k,l,m,t)=\f{1}{2\om_j}\big(e^{-i\om_j
t}v^{(j)+}(k,l,m,t)-e^{i\om_j t}v^{(j)-}(k,l,m,t)\big),
 \end{align}
\begin{flushright}
$k,l\in\mathbb{R}, m\in\mathbb{C}, t\geq0, j=1,2.$\end{flushright}
\end{subequations}
Similarly, let $g^{(j)\pm}$, $f^{(j)\pm}$, $e^{(j)\pm}$, $G^{(j)}$,
$F^{(j)},~E^{(j)},$ $j=1,2,$ denote the following \textit{known} functions:
\begin{subequations}\label{BigG}
 \begin{align}
&g^{(j)\pm}(k,l,m,t)=\int_0^t e^{\pm i\om_j s}\q{g}_1(k,l,s)ds,\\
&f^{(j)\pm}(k,l,m,t)=\int_0^t e^{\pm i\om_j s}\q{g}_2(k,l,s)ds,\\
&e^{(j)\pm}(k,l,m,t)=\int_0^t e^{\pm i\om_j s}\q{g}_3(k,l,s)ds,\\
&G^{(j)}(k,l,m,t)=\f{1}{2\om_j}\big(e^{-i\om_j t}g^{(j)+}(k,l,t)-e^{i\om_j t}g^{(j)-}(k,l,m,t)\big),\\
&F^{(j)}(k,l,m,t)=\f{1}{2\om_j}\big(e^{-i\om_j t}f^{(j)+}(k,l,t)-e^{i\om_j t}f^{(j)-}(k,l,m,t)\big),\\
&E^{(j)}(k,l,m,t)=\f{1}{2\om_j}\big(e^{-i\om_j t}e^{(j)+}(k,l,t)-e^{i\om_j t}e^{(j)-}(k,l,m,t)\big),
\end{align}
\begin{flushright}$k\in\mathbb{R}, l\in\mathbb{C}, t\geq0, j=1,2.$\end{flushright}
\end{subequations}

Using the above notations, equations \eqref{global_r} become
\begin{subequations}\label{global_relations}
 \begin{align}
\label{1}\begin{split}
&k\hat{u}+l\hat{v}+m\hat{w}=2\mu km U^{(1)}+2\mu m l V^{(1)}+(\la(k^2+l^2)\\
&~~~~~~~~~~~~~~~~~~~~~~~~~~~+m^2(\la+2\mu))W^{(1)}+N_P,
\end{split}\\
\label{2}
& m\hat{v}-l\hat{w}=-\mu kl U^{(2)}+\mu(m^2-l^2) V^{(2)}-2\mu ml W^{(2)}+N_Q,\\
\label{3}
& k\hat{w}-m\hat{u}=\mu (k^2-m^2) U^{(2)}+\mu kl V^{(2)}+2\mu km W^{(2)}+N_R,
\end{align}
\end{subequations}
where the \textit{known} functions $N_P(k,l,m,t)$, $N_Q(k,l,m,t)$ and $N_R(k,l,m,t)$ are defined as follows:
\begin{subequations}\label{NPNQ}
 \begin{align}
\begin{split}
&N_P(k,l,m,t)= -i\mu k G^{(1)}(k,l,m,t) -i\mu l F^{(1)}(k,l,m,t)-im(\la+2\mu) E^{(1)}(k,l,m,t)\\
&~~+\Big(\f{1}{2}P_0(k,l,m)+\f{i}{2\om_1}P_1(k,l,m)\Big)e^{-i\om_1t}+\Big(\f{1}{2}P_0(k,l,m)-\f{i }{2\om_1}P_1(k,l,m)\Big)e^{i\om_1t},
\end{split}\\
\begin{split}
&N_Q(k,l,m,t)=-i\mu m F^{(2)}(k,l,m,t)- il(\la+2\mu) E^{(2)}(k,l,m,t)\\
&~~+\Big(\f{1}{2}Q_0(k,l,m)+\f{i}{2\om_2}Q_1(k,l,m)\Big)e^{-i\om_2t}+\Big(\f{1}{2}Q_0(k,l,m)-\f{i}{2\om_2}Q_1(k,l,m)\Big)e^{i\om_2t}.
\end{split}\\
\begin{split}
&N_R(k,l,m,t)=i\mu m G^{(2)}(k,l,m,t)-ik(\la+2\mu) E^{(2)}(k,l,m,t)\\
&~~+\Big(\f{1}{2}R_0(k,l,m)+\f{i}{2\om_2}R_1(k,l,m)\Big)e^{-i\om_2t}+\Big(\f{1}{2}R_0(k,l,m)-\f{i}{2\om_2}R_1(k,l,m)\Big)e^{i\om_2t}.
\end{split}
\end{align}
\end{subequations}
We will refer to equations \eqref{global_relations} as the \textit{global relations}.
The global relations express the three dimensional Fourier
transforms of the solution $(u,v,w)$ in terms of the given initial and
boundary data, as well as in terms of the transforms
$U^{(1)},U^{(2)},V^{(1)},V^{(2)},$ $W^{(1)},W^{(2)}$ of the unknown boundary values. The important
observation is that equations \eqref{global_relations} are valid for all
values of $m$ in the lower half complex $m-$plane. It turns out that
using this fact it will be possible to eliminate the unknown
transforms.

\section{An Integral Representation Involving the Unknown Transforms}
We first observe that for any fixed $k,l\in\mathbb{R}$ and any fixed
$t$, $0\leq t< T$, $T>0$, the functions $U^{(j)}$, $V^{(j)}$, $W^{(j)}$,
$G^{(j)},~F^{(j)},$ $E^{(j)}$, $j\in\{1,2\},$ are analytic in the complex
$m-$plane.

Indeed, these functions are all single-valued, thus it only remains
to establish the analyticity in the neighbourhood $m=\pm i\sqrt{k^2+l^2}$. The
definition of $W^{(1)}$, i.e. equation \eqref{defW1}, implies that
this function posses the following expansion for $m$ near $\pm i\sqrt{k^2+l^2}$:
\begin{equation}\label{expansion}
W^{(1)}(k,l,m,t)=i\sum_{n=1}^{\infty}\f{\om_1^{2n-2}}{(2n-1)!}\int_0^t(s-t)^{2n-2}\q{w}(k,l,s)ds.
\end{equation}
Similar expansions are also valid for $U^{(1)},U^{(2)},V^{(1)},V^{(2)},$ $W^{(2)},$
$G^{(1)},G^{(2)},$ $F^{(1)},F^{(2)}$, $E^{(1)},E^{(2)}.$

%The following symmetry relations are valid:
%$$u^{(j)\pm}(k,l,m,t)=u^{(j)\mp}(k,l,-m,t),~v^{(j)\pm}(k,l,m,t)=v^{(j)\mp}(k,l,-m,t),$$
%$$w^{(j)\pm}(k,l,m,t)=w^{(j)\mp}(k,l,-m,t),$$
%$$U^{(j)}(k,l,m,t)=U^{(j)}(k,l,-m,t),~~V^{(j)}(k,l,m,t)=V^{(j)}(k,l,-m,t),$$
%$$W^{(j)}(k,l,m,t)=W^{(j)}(k,l,-m,t),~~~j=1,2.$$
Solving equations \eqref{global_relations} we find
\begin{subequations}
 \label{gwh}
 \begin{align}
\label{gwh1}
\begin{split}
&\hat{u}=\f{1}{m(k^2+l^2+m^2)}\big[2\mu k^2m^2 U^{(1)}+2\mu m^2 kl V^{(1)}+[\la km(k^2+l^2)+km^3(\la+2\mu)]W^{(1)}\\
&-\mu m^2 (k^2-l^2-m^2) U^{(2)}-2\mu m^2 kl V^{(2)}-2\mu k m^3 W^{(2)}+km N_P-kl N_Q-(m^2+l^2)N_R\big],
\end{split}\\
\label{gwh2}
\begin{split}
&\hat{v}=\f{1}{m(k^2+l^2+m^2)}\big[2\mu kl m^2 U^{(1)}+ 2\mu m^2 l^2 V^{(1)} + [\la lm (k^2+l^2)+lm^3(\la+2\mu)]W^{(1)}\\
&-2\mu m^2 kl U^{(2)} + \mu m^2 (m^2+k^2-l^2) V^{(2)}-2\mu l m^3 W^{(2)}+lm N_P+(k^2+m^2)N_Q+lkN_R\big],
\end{split}\\
\label{gwh3}
\begin{split}
&\hat{w}=\f{1}{k^2+l^2+m^2}\big[2\mu k m^2 U^{(1)}+ 2\mu l m^2 V^{(1)} + [\la k^2 m+\la l^2 m+m^3(\la+2\mu)]W^{(1)}\\
&+\mu k (l^2+k^2-m^2)U^{(2)} + \mu l (k^2+l^2-m^2) V^{(2)}+2\mu m (k^2+l^2) W^{(2)}+m N_P-l N_Q +k N_R\big].
\end{split}
%\begin{split}
%\end{split}
 \end{align}
\end{subequations}

We observe that the following important transformation is valid in
the complex $m-$plane:
\begin{equation}\label{superb_1}
\begin{split}
m\rightarrow -m~~&maps~~\om_1~~ to~~ -\om_1,~~~ \om_2~~ to~~ -\om_2,~
and~~\\
&leaves~~ U^{(j)},V^{(j)},W^{(j)},~ j=1,2~~ invariant.
\end{split}
\end{equation}
Employing the transformation $m\rightarrow-m$ in equations
\eqref{gwh} and then adding the resulting
equations to \eqref{gwh}, we obtain the following
equations:
\begin{subequations}
 \begin{align}
\label{u}
\begin{split}
&\hat{u}(k,l,m,t)+\hat{u}(k,l,-m,t)=\\
&\f{1}{m(k^2+l^2+m^2)}\big[2km[\la (k^2+l^2)+m^2(\la+2\mu)]W^{(1)}-4\mu k m^3 W^{(2)}+\\
&km[N_P(k,l,m,t)+N_P(k,l,-m,t)]-kl[N_Q(k,l,m,t)-N_Q(k,l,-m,t)]-\\
&(m^2+l^2)[N_R(k,l,m,t)-N_R(k,l,-m,t)]\big],
\end{split}\\
\label{v}
\begin{split}
&\hat{v}(k,l,m,t)+\hat{v}(k,l,-m,t)=\\
&\f{1}{m(k^2+l^2+m^2)}\big[ 2lm [\la (k^2+l^2)+m^2(\la+2\mu)]W^{(1)}-4\mu l m^3 W^{(2)}+\\
&lm [N_P(k,l,m,t)+N_P(k,l,-m,t)]+(k^2+m^2)[N_Q(k,l,m,t)-\\
&N_Q(k,l,-m,t)]+lk[N_R(k,l,m,t)-N_R(k,l,-m,t)]\big],
\end{split}\\
\begin{split}
&\hat{w}(k,l,m,t)+\hat{w}(k,l,-m,t)=\\
&\f{1}{k^2+l^2+m^2}\big[4\mu k m^2 U^{(1)}+ 4\mu l m^2 V^{(1)}+2\mu k (l^2+k^2-m^2)U^{(2)} +\\
& 2\mu l (k^2+l^2-m^2) V^{(2)}+m[N_P(k,l,m,t)-N_P(k,l,-m,t)]-\\
&l[N_Q(k,l,m,t)+N_Q(k,l,-m,t)]+k[N_R(k,l,m,t)+N_R(k,l,-m,t)]\big].
\end{split}
 \end{align}
\end{subequations}
Applying the inverse Fourier transform formula to these equations,
we obtain
\begin{subequations}\label{first}
 \begin{align}
\label{first_u}
\begin{split}
&u(x,y,z,t)=\f{1}{8\pi^3}\int_{-\infty}^{\infty}dk \int_{-\infty}^{\infty} dl \int_{-\infty}^{\infty} dm~e^{ikx+ily+imz}\hat{u}(k,l,m,t)\\
&=\f{1}{8\pi^3}\int_{-\infty}^{\infty}\!\!dk
\int_{-\infty}^{\infty}\!\! dl \int_{-\infty}^{\infty}\!\! dm
\f{e^{ikx+ily+imz}}{k^2+l^2+m^2}\\
& \Big\{ 2k(\la (k^2+l^2)+m^2(\la+2\mu))W^{(1)}-4\mu k m^2 W^{(2)}+\\
&k[N_P(k,l,m,t)+N_P(k,l,-m,t)]-\f{1}{m}kl[N_Q(k,l,m,t)-N_Q(k,l,-m,t)]-\\
&\f{1}{m}(m^2+l^2)[N_R(k,l,m,t)-N_R(k,l,-m,t)]\Big\},
\end{split}\\
\label{first_v}
\begin{split}
&v(x,y,z,t)=\f{1}{8\pi^3}\int_{-\infty}^{\infty}dk \int_{-\infty}^{\infty} dl \int_{-\infty}^{\infty} dm~e^{ikx+ily+imz}\hat{v}(k,l,m,t)\\
&=\f{1}{8\pi^3}\int_{-\infty}^{\infty}\!\!dk
\int_{-\infty}^{\infty}\!\! dl \int_{-\infty}^{\infty}\!\! dm
\f{e^{ikx+ily+imz}}{k^2+l^2+m^2}\\
& \Big\{ 2l (\la (k^2+l^2)+m^2(\la+2\mu))W^{(1)}-4\mu l m^2 W^{(2)}+\\
& l [N_P(k,l,m,t)+N_P(k,l,-m,t)]+\f{1}{m}(k^2+m^2)[N_Q(k,l,m,t)-\\
& N_Q(k,l,-m,t)]+\f{1}{m}lk[N_R(k,l,m,t)-N_R(k,l,-m,t)]\Big\},
\end{split}\\
\begin{split}
&w(x,y,z,t)=\f{1}{8\pi^3}\int_{-\infty}^{\infty}dk \int_{-\infty}^{\infty} dl \int_{-\infty}^{\infty} dm~e^{ikx+ily+imz}\hat{w}(k,l,m,t)\\
&=\f{1}{8\pi^3}\int_{-\infty}^{\infty}\!\!dk
\int_{-\infty}^{\infty}\!\! dl \int_{-\infty}^{\infty}\!\! dm
\f{e^{ikx+ily+imz}}{k^2+l^2+m^2}\\
& \Big\{4\mu k m^2 U^{(1)}+ 4\mu l m^2 V^{(1)}+2\mu k (l^2+k^2-m^2)U^{(2)} +\\
& 2\mu l (k^2+l^2-m^2) V^{(2)}+m[N_P(k,l,m,t)-N_P(k,l,-m,t)]-\\
& l[N_Q(k,l,m,t)+N_Q(k,l,-m,t)]+k[N_R(k,l,m,t)+N_R(k,l,-m,t)]\Big\},
\end{split}
 \end{align}
\begin{flushright}$-\infty<x<\infty,~-\infty<y<\infty,~0<z<\infty,~t>0.$\end{flushright}
\end{subequations}

Denote by $H_j(k,l,m,t)$, $j=1,2,3,$ the following functions
appearing in equations \eqref{first}:
\begin{subequations}
 \begin{align}
&H_1(k,l,m,t)=\f{1}{k^2+l^2+m^2}\Big\{ 2k[\la (k^2+l^2)+m^2(\la+2\mu)]W^{(1)}-4\mu k m^2 W^{(2)}\Big\},\\
&H_2(k,l,m,t)=\f{1}{k^2+l^2+m^2}\Big\{ 2l[\la (k^2+l^2)+m^2(\la+2\mu)]W^{(1)}-4\mu l m^2 W^{(2)}\Big\},\\
\begin{split}
&H_3(k,l,m,t)=\f{1}{k^2+l^2+m^2}\Big\{ 4\mu k m^2 U^{(1)}+ 4\mu l m^2 V^{(1)}\\
&~~~~~~~~~~~~~~~~~~~~~~~~~~+2\mu k (l^2+k^2-m^2)U^{(2)} +2\mu l (k^2+l^2-m^2) V^{(2)}\Big\}.
\end{split}
 \end{align}
\end{subequations}
We observe that for any fixed $k,l\in\mathbb{R}$ and fixed $t$, $0\leq
t<T$, $T>0$, the above functions are analytic in the entire complex
$m-$plane. Indeed, equation \eqref{expansion} and the analogous
equation for $W^{(2)}$, imply that in the neighbourhood of $m=\pm
i\sqrt{k^2+l^2}$, the following expansion is valid:
\begin{equation}
H_1(k,l,m,t)=2i\la k \int_0^t\q{u}(k,l,s)ds+o(\om_1^2).
\end{equation}
Similarly,
\begin{equation}
H_2(k,l,m,t)=2i\la l \int_0^t\q{v}(k,l,s)ds+o(\om_1^2),
\end{equation}
\begin{equation}
H_3(k,l,m,t)=2i\mu k \int_0^t\q{u}(k,l,s)ds+2i\mu l \int_0^t\q{v}(k,l,s)ds+o(\om_1^2).
\end{equation}

The restriction $z>0$, as well as the analyticity of the functions
$H_j(k,l,m,t)$ $j=1,2,3,$ allow us to deform the contour of
integration from the real axis to a contour $\ga_{k,l}$ in the upper
half $m-$plane (the particular choice of $\ga_{k,l}$ will be determined
in the next section):
\begin{subequations}\label{before_final}
 \begin{align}
\label{first_uf}
\begin{split}
&u(x,y,z,t)=\f{1}{8\pi^3}\int_{-\infty}^{\infty}dk \int_{-\infty}^{\infty} dl \int_{\gamma_{k,l}} dm~\f{e^{ikx+ily+imz}}{k^2+l^2+m^2}\\
& \Big\{ 2k[\la (k^2+l^2)+m^2(\la+2\mu)]W^{(1)}-4\mu k m^2 W^{(2)}\Big\}\\
& +\f{1}{8\pi^3}\int_{-\infty}^{\infty}\!\!dk \int_{-\infty}^{\infty}\!\! dl \int_{-\infty}^{\infty}\!\! dm~ \f{e^{ikx+ily+imz}}{k^2+l^2+m^2}\\
& \Big\{k[N_P(k,l,m,t)+N_P(k,l,-m,t)]-\f{1}{m}kl[N_Q(k,l,m,t)-N_Q(k,l,-m,t)]-\\
& \f{1}{m}(m^2+l^2)[N_R(k,l,m,t)-N_R(k,l,-m,t)]\Big\},
\end{split}\\
\label{first_vf}
\begin{split}
& v(x,y,z,t)=\f{1}{8\pi^3}\int_{-\infty}^{\infty}dk \int_{-\infty}^{\infty} dl \int_{\gamma_{k,l}} dm~\f{e^{ikx+ily+imz}}{k^2+l^2+m^2}\\
& \Big\{ 2l [\la (k^2+l^2)+m^2(\la+2\mu)]W^{(1)}-4\mu l m^2 W^{(2)}\Big\}\\
& +\f{1}{8\pi^3}\int_{-\infty}^{\infty}\!\!dk \int_{-\infty}^{\infty}\!\! dl \int_{-\infty}^{\infty}\!\! dm~ \f{e^{ikx+ily+imz}}{k^2+l^2+m^2}\\
& \Big\{l [N_P(k,l,m,t)+N_P(k,l,-m,t)]+\f{1}{m}(k^2+m^2)[N_Q(k,l,m,t)-\\
& N_Q(k,l,-m,t)]+\f{1}{m}lk[N_R(k,l,m,t)-N_R(k,l,-m,t)]\Big\},
\end{split}\\
\label{first_wf}
\begin{split}
& w(x,y,z,t)=\f{1}{8\pi^3}\int_{-\infty}^{\infty}dk \int_{-\infty}^{\infty} dl \int_{\gamma_{k,l}} dm~\f{e^{ikx+ily+imz}}{k^2+l^2+m^2}\\
& \Big\{ 4\mu k m^2 U^{(1)}+ 4\mu l m^2 V^{(1)}+2\mu k (l^2+k^2-m^2)U^{(2)}+ 2\mu l (k^2+l^2-m^2) V^{(2)}\Big\}\\
& +\f{1}{8\pi^3}\int_{-\infty}^{\infty}\!\!dk \int_{-\infty}^{\infty}\!\! dl \int_{-\infty}^{\infty}\!\! dm~ \f{e^{ikx+ily+imz}}{k^2+l^2+m^2}\\
& \Big\{m[N_P(k,l,m,t)-N_P(k,l,-m,t)]-\\
& l[N_Q(k,l,m,t)+N_Q(k,l,-m,t)]+k[N_R(k,l,m,t)+N_R(k,l,-m,t)]\Big\}.
\end{split}
\end{align}
\end{subequations}

\section{The Elimination of the Transforms of the Unknown Boundary Values}
Let
\begin{equation}\label{m21}
m_{21}=-m\Big(\f{\la+2\mu}{\mu}+\f{\la+\mu}{\mu}\f{k^2+l^2}{m^2}\Big)^{\f{1}{2}}
\end{equation}
and
\begin{equation}\label{m12}
m_{12}=-m\Big(\f{\mu}{\la+2\mu}-\f{\la+\mu}{\la+2\mu}\f{k^2+l^2}{m^2}\Big)^{\f{1}{2}}.
\end{equation}

\begin{figure}[ht]
\begin{center} \small
$$
\begin{array}{cc}
\resizebox{5.5cm}{!}{\input{l12.pstex_t}} & \resizebox{5.5cm}{!}{\input{l21.pstex_t}}\\
\mbox{a. The branch cut for }m_{12}. & \mbox{b. The branch cut for } m_{21}.
\end{array}
$$
\end{center}
\caption{Branch cuts for $m_{12}$ and $m_{21}$ in complex
$m-$plane.} \label{m12l21}
\end{figure}
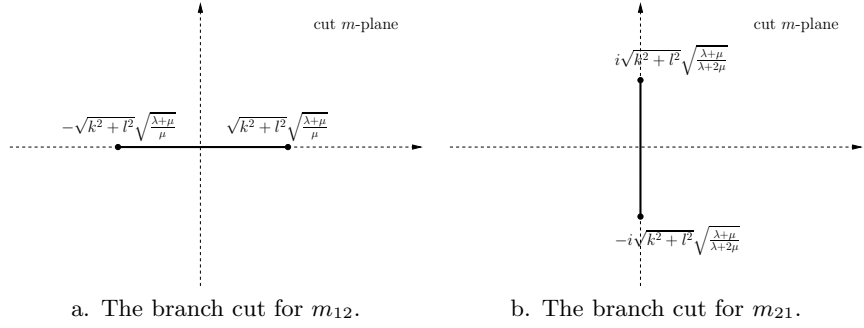
The function $m_{12}$ has two branch points $\pm
\sqrt{k^2+l^2}\sqrt{\f{\la+\mu}{\mu}}$, which we connect by a horizontal branch
cut; the function $m_{21}$ has two branch points $\pm
i\sqrt{k^2+l^2}\sqrt{\f{\la+\mu}{\la+2\mu}}$, which we connect by a vertical
branch cut, see Fig.\ref{m12l21}. We fix the branches by the
requirements that,
$$m_{21}\sim -m\sqrt{\f{\la+2\mu}{\mu}},~~m_{12}\sim-m\sqrt{\f{\mu}{\la+2\mu}},~~ as~m\rightarrow\infty.$$

The following transformations are valid in the cut $m$-plane:
\begin{subequations}\label{superb_2}
\begin{align}
&m\rightarrow m_{21}~~ maps~~~ \om_2~~ to~~ -\om_1,~~~ U^{(2)}~~
to~~ U^{(1)},~~~ V^{(2)}~~ to~~V^{(1)},~~~W^{(2)}~~to~~W^{(1)};\\
&m\rightarrow m_{12}~~ maps~~~ \om_1~~ to~~ -\om_2,~~~ U^{(1)}~~
to~~ U^{(2)},~~~ V^{(1)}~~ to~~ V^{(2)},~~~ W^{(1)}~~ to~~ W^{(2)}.
\end{align}
\end{subequations}
%\paragraph{(ii).} Elimination of the unknowns.

Using in equations \eqref{1},\eqref{2},\eqref{3} the transformations
$m\rightarrow m_{12}$, $m\rightarrow-m$, $m\rightarrow-m$, respectively, and then
combining the three resulting equations, we obtain the following equations:
% \begin{align}
%\label{1}
%&k\hat{u}+l\hat{v}+m\hat{w}=2\mu km U^{(1)}+2\mu m l V^{(1)}+(\la(k^2+l^2)+m^2(\la+2\mu))W^{(1)}+N_P,\\
%\label{2}
%&m\hat{v}-l\hat{w}=-\mu kl U^{(2)}+\mu(m^2-l^2) V^{(2)}-2\mu ml W^{(2)}+N_Q,\\
%\label{3}
%&k\hat{w}-m\hat{u}=\mu (k^2-m^2) U^{(2)}+\mu kl V^{(2)}+2\mu km W^{(2)}+N_R,
% \end{align}
\begin{subequations}\label{C1}
 \begin{align}
 \begin{split}
&k\hat{u}(k,l,m_{12},t)+l\hat{v}(k,l,m_{12},t)+m_{12}\hat{w}(k,l,m_{12},t)=2\mu km_{12} U^{(2)}(k,l,m,t)+\\
&2\mu m_{12} l V^{(2)}(k,l,m,t)+(\la(k^2+l^2)+m_{12}^2(\la+2\mu))W^{(2)}(k,l,m,t)+N_P(k,l,m_{12},t)\\
\end{split}\\
 \begin{split}
&-m\hat{v}(k,l,-m,t)-l\hat{w}(k,l,-m,t)=-\mu kl U^{(2)}(k,l,m,t)+\mu(m^2-l^2) V^{(2)}(k,l,m,t)+\\
&2\mu ml W^{(2)}(k,l,m,t)+N_Q(k,l,-m,t),
\end{split}\\
\begin{split}
&k\hat{w}(k,l,-m,t)+m\hat{u}(k,l,-m,t)=\mu (k^2-m^2) U^{(2)}(k,l,-m,t)+\mu kl V^{(2)}(k,l,-m,t)-\\
&2\mu km W^{(2)}(k,l,-m,t)+N_R(k,l,-m,t).
\end{split}
 \end{align}
\end{subequations}

Similarly, using in equations \eqref{1}, \eqref{2}, \eqref{3} the
transformations $m\rightarrow -m$, $m\rightarrow m_{21}$, $m\rightarrow m_{21}$
respectively, and then combining the three resulting equations, we
obtain the following equations:
\begin{subequations}\label{C2}
 \begin{align}
 \begin{split}
  &k\hat{u}(k,l,-m,t)+l\hat{v}(k,l,-m,t)-m\hat{w}(k,l,-m,t)=-2\mu km U^{(1)}(k,l,-m,t)-\\
  &2\mu m l V^{(1)}(k,l,-m,t)+(\la(k^2+l^2)+m^2(\la+2\mu))W^{(1)}(k,l,-m,t)+N_P(k,l,-m,t),
 \end{split}\\
 \begin{split}
  &m_{21}\hat{v}(k,l,m_{21},t)-l\hat{w}(k,l,m_{21},t)=-\mu kl U^{(1)}(k,l,m,t)+\\
  &\mu(m_{21}^2-l^2) V^{(1)}(k,l,m,t)-2\mu l m_{21} W^{(1)}(k,l,m,t)+N_Q(k,l,m_{21},t),
 \end{split}\\
  \begin{split}
   &k\hat{w}(k,l,m_{21},t)-m_{21}\hat{u}(k,l,m_{21},t)=\mu (k^2-m_{21}^2) U^{(1)}(k,l,m,t)+\\
   &\mu kl V^{(1)}(k,l,m,t)+2\mu km_{21} W^{(1)}(k,l,m,t)+N_R(k,l,m_{21},t).
  \end{split}
 \end{align}
\end{subequations}

Let
\begin{subequations}\label{Cs}
\begin{align}
&C=\left(
\begin{array}{llcl}
 2\mu km_{12}               & 2\mu m_{12} l    & \la(k^2+l^2)+m_{12}^2(\la+2\mu)\\
 -\mu kl                    & \mu(m^2-l^2)     & 2\mu ml \\
 -\mu (m^2-k^2)              & \mu kl           & -2\mu km \\
\end{array}\right),\\
&D=\left(
\begin{array}{llcl}
 -2\mu km                   & -2\mu m l              & \la(k^2+l^2)+m^2(\la+2\mu)\\
 -\mu kl                    & \mu(m_{21}^2-l^2)     & -2\mu lm_{21} \\
 \mu (k^2-m_{21}^2)         & \mu kl                & 2\mu km_{21} \\
\end{array}\right),\\
&\Delta_1=det(C),~~~\Delta_2=det(D).
\end{align}
\end{subequations}
Simplifying the expressions for $\Delta_1$ and $\Delta_2$ we find
\begin{subequations}\label{Ts}
\begin{align}
\Delta_1&=\mu^3 m^2\big[(k^2+l^2-m^2)^2-4(k^2+l^2)m m_{12}\big],\\
\Delta_2&=\mu m_{21}^2\big[(\la(k^2+l^2)+m^2(\la+2\mu))^2-4\mu^2 (k^2+l^2) mm_{21}\big].
\end{align}
\end{subequations}
Equations \eqref{C1} imply
\begin{equation}\label{X}
\left(
\begin{array}{llcl}
 U^{(2)}(k,l,m,t)\\
 V^{(2)}(k,l,m,t)\\
 W^{(2)}(k,l,m,t)\\
\end{array}\right)=C^{-1}
\left(
\begin{array}{llcl}
k\hat{u}(k,l,m_{12},t)+l\hat{v}(k,l,m_{12},t)+m_{12}\hat{w}(k,l,m_{12},t)-N_P(k,l,m_{12},t)\\
-m\hat{v}(k,l,-m,t)-l\hat{w}(k,l,-m,t)-N_Q(k,l,-m,t)\\
k\hat{w}(k,l,-m,t)+m\hat{u}(k,l,-m,t)-N_R(k,l,-m,t)\\
\end{array}\right),
\end{equation}
Equations \eqref{C2} imply
\begin{equation}\label{Y}
\left(
\begin{array}{llcl}
 U^{(1)}(k,l,m,t)\\
 V^{(1)}(k,l,m,t)\\
 W^{(1)}(k,l,m,t)\\
\end{array}\right)=D^{-1}
\left(
\begin{array}{llcl}
k\hat{u}(k,l,-m,t)+l\hat{v}(k,l,-m,t)-m\hat{w}(k,l,-m,t)-N_P(k,l,-m,t)\\
m_{21}\hat{v}(k,l,m_{21},t)-l\hat{w}(k,l,m_{21},t)-N_Q(k,l,m_{21},t)\\
k\hat{w}(k,l,m_{21},t)-m_{21}\hat{u}(k,l,m_{21},t)-N_R(k,l,m_{21},t)\\
\end{array}\right),
\end{equation}
where $C^{-1}$ and $D^{-1}$ are the inverse matrices of $C$ and $D$ respectively.

We fix the choice of the contour $\ga_{k,l}$ by requiring that every
term in the RHS of \eqref{X} and \eqref{Y} does not have a pole or a
branch point above this contour, see Fig.\ref{gamma_{k,l}}.
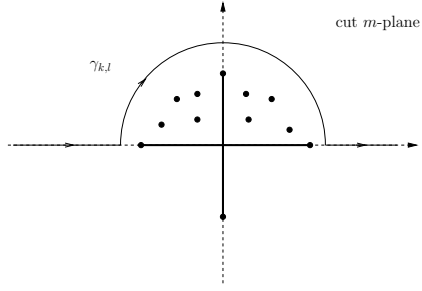
\begin{figure}[ht]
\begin{center}
\resizebox{5.5cm}{!}{\input{Newversion_1.pstex_t}}
\end{center}
\caption{$\gamma_{k,l}$, the deformed path of integration.} \label{gamma_{k,l}}
\end{figure}

Regarding the zeros of $\Delta_j$, $j=1,2$, we note that they are of the form
$m=\al \sqrt{k^2+l^2}$, for some constant $\al\in \mathbb{C}$. For example, if $\la=2\mu$, we find
that the zeros of $\Delta_1$ are
$$m=0,~m=0,~m\approx (-1.624\pm 0.126 i)\sqrt{k^2+l^2},~m\approx \pm 0.357 i\sqrt{k^2+l^2},$$
$$m\approx \pm 1.056 i\sqrt{k^2+l^2},~m \approx (1.624\pm 0.126i)\sqrt{k^2+l^2};$$
whereas the zeros of $\Delta_2$ are
$$m\approx \pm 0.866 i \sqrt{k^2+l^2},~m\approx (-0.295\pm 0.442 i)\sqrt{k^2+l^2},~m\approx \pm 0.885 i\sqrt{k^2+l^2},$$
$$m\approx \pm i\sqrt{k^2+l^2},~m \approx (0.295\pm 0.442i)\sqrt{k^2+l^2}.$$

Substituting equations \eqref{X} and \eqref{Y} in equations
\eqref{before_final},
%although there may be poles from the zeros of
%$\Delta_1,\Delta_2$; there may also be branch points from
%$\Delta_1,\Delta_2$ or $C_2$ and $D_3$.
and using Jordan's lemma in the complex $m-$plane above the contour
$\ga_k$, it follows that $\hat{u}(k,l,-m,t)$, $\hat{u}(k,l,m_{12},t)$,
$\hat{u}(k,l,m_{21},t)$, $\hat{v}(k,l,-m,t)$, $\hat{v}(k,l,m_{12},t)$,
$\hat{v}(k,l,m_{21},t)$, $\hat{w}(k,l,-m,t)$, $\hat{w}(k,l,m_{12},t)$,
$\hat{w}(k,l,m_{21},t)$ yield a zero contribution. Let
\begin{subequations}\label{notation}\begin{align}
\left(
\begin{array}{llcl}
 U_2(k,l,m,t)\\
 V_2(k,l,m,t)\\
 W_2(k,l,m,t)\\
\end{array}\right)=-C^{-1}
\left(
\begin{array}{llcl}
N_P(k,l,m_{12},t)\\
N_Q(k,l,-m,t)\\
N_R(k,l,-m,t)\\
\end{array}\right),\\
\left(
\begin{array}{llcl}
 U_1(k,l,m,t)\\
 V_1(k,l,m,t)\\
 W_1(k,l,m,t)\\
\end{array}\right)=-D^{-1}
\left(
\begin{array}{llcl}
N_P(k,l,-m,t)\\
N_Q(k,l,m_{21},t)\\
N_R(k,l,m_{21},t)\\
\end{array}\right),
\end{align}
\end{subequations}
where $N_P,N_Q,N_R$ are known functions defined in \eqref{NPNQ}, and $(m_{12},m_{21})$ are defined by equations \eqref{m12} and \eqref{m21} respectively.
Equations \eqref{before_final} become:
\begin{subequations}\label{final}
 \begin{align}
\label{final_u}
\begin{split}
&u(x,y,z,t)=\\
&\f{1}{8\pi^3}\int_{-\infty}^{\infty}dk \int_{-\infty}^{\infty} dl \int_{\gamma_{k,l}} dm~\f{e^{ikx+ily+imz}}{k^2+l^2+m^2}\\
&\Big\{ 2k[\la (k^2+l^2)+m^2(\la+2\mu)]W_1-4\mu k m^2 W_2\Big\}\\
& +\f{1}{8\pi^3}\int_{-\infty}^{\infty}\!\!dk \int_{-\infty}^{\infty}\!\! dl \int_{-\infty}^{\infty}\!\! dm \f{e^{ikx+ily+imz}}{k^2+l^2+m^2}\\
& \Big\{k[N_P(k,l,m,t)+N_P(k,l,-m,t)]-\f{1}{m}kl[N_Q(k,l,m,t)-N_Q(k,l,-m,t)]-\\
& \f{1}{m}(m^2+l^2)[N_R(k,l,m,t)-N_R(k,l,-m,t)]\Big\},
\end{split}\\
\label{final_v}
\begin{split}
& v(x,y,z,t)=\\
&\f{1}{8\pi^3}\int_{-\infty}^{\infty}dk \int_{-\infty}^{\infty} dl \int_{\gamma_{k,l}} dm~\f{e^{ikx+ily+imz}}{k^2+l^2+m^2}\\
& \Big\{ 2l [\la (k^2+l^2)+m^2(\la+2\mu)]W_1-4\mu l m^2 W_2\Big\}\\
& +\f{1}{8\pi^3}\int_{-\infty}^{\infty}\!\!dk \int_{-\infty}^{\infty}\!\! dl \int_{-\infty}^{\infty}\!\! dm \f{e^{ikx+ily+imz}}{k^2+l^2+m^2}\\
& \Big\{l [N_P(k,l,m,t)+N_P(k,l,-m,t)]+\f{1}{m}(k^2+m^2)[N_Q(k,l,m,t)-\\
& N_Q(k,l,-m,t)]+\f{1}{m}lk[N_R(k,l,m,t)-N_R(k,l,-m,t)]\Big\},
\end{split}\\
\label{final_w}
\begin{split}
& w(x,y,z,t)=\\
&\f{1}{8\pi^3}\int_{-\infty}^{\infty}dk \int_{-\infty}^{\infty} dl \int_{\gamma_{k,l}} dm~\f{e^{ikx+ily+imz}}{k^2+l^2+m^2}\\
& \Big\{ 4\mu k m^2 U_1+ 4\mu l m^2 V_1+2\mu k (l^2+k^2-m^2)U_2+ 2\mu l (k^2+l^2-m^2) V_2\Big\}\\
& +\f{1}{8\pi^3}\int_{-\infty}^{\infty}\!\!dk \int_{-\infty}^{\infty}\!\! dl \int_{-\infty}^{\infty}\!\! dm \f{e^{ikx+ily+imz}}{k^2+l^2+m^2}\\
& \Big\{m[N_P(k,l,m,t)-N_P(k,l,-m,t)]-\\
& l[N_Q(k,l,m,t)+N_Q(k,l,-m,t)]+k[N_R(k,l,m,t)+N_R(k,l,-m,t)]\Big\}.
\end{split}
\end{align}
\end{subequations}
where $C,D$ are defined in \eqref{Cs}, the known functions $U_j,V_j,W_j$, $j=1,2$, are defined in \eqref{notation}, and the known functions $N_P,N_Q,N_R$ are defined in \eqref{NPNQ}.

We summarize the above result in the following proposition:
\begin{proposition}
Let $(u,v,w)$ satisfy the Lam\'e-Navier equations \eqref{3D} in the half space
$$-\infty<x<\infty,~-\infty<y<\infty,~0<z<\infty,~t>0,$$
with the initial conditions \eqref{ic} and the stress boundary conditions \eqref{3D_l_bv}. Assume that the given
functions $$\{u_j(x,y,z),v_j(x,y,z),w_j(x,y,z)\}_{j=0,1},~~\{g_j(x,y,t)\}_{j=1,2,3}$$ have sufficient smoothness and decay.
A solution of the above initial-boundary value problem, which decays for large $(x,y,z)$, is given by equations \eqref{final}, where:
\begin{enumerate}
\item[a)] The known functions $(N_P,N_Q,N_R)$ are defined in \eqref{NPNQ} in terms of the transforms of the initial and boundary data (see equations \eqref{dispersion}-\eqref{BigG}).
\item[b)] The known functions $(U_j,V_j,W_j)$ are defined in \eqref{notation} in terms of $(N_P,$ $N_Q,N_R)$ and of the matrices $C$ and $D$ given by \eqref{Cs}.
\item[c)] The contours $\gamma_{k,l}$ depicted in Fig.\ref{gamma_{k,l}}, are deformations of the real axis and determined by the requirement that the zeros of $\Delta_j$, $j=1,2$ and the branch points of $m_{21}$ and $m_{12}$ are all below $\gamma_{k,l}$.
\end{enumerate}
\end{proposition}

\section{The Normal Point Load with Homogeneous Initial Conditions}
Consider a normal point load suddenly applied to an isotropic elastic half space body.
In this case,
\begin{equation}\label{Lamb_normal_point}
u_0=u_1=v_0=v_1=w_0=w_1=0; g_1=0,~g_2=0,~g_3=\sigma_0\delta(x,y)h(t)/(\la+\mu),
\end{equation}
where $\sigma_0$ is a constant and $h(t)$ is the Heaviside
function defined by $h(t)=0,~t\leq 0;$ $h(t)=1,~t>0$.
%For a moving normal line load with a constant speed $C$,
%$g_1=0,~g_2=\sigma_0\delta(x-Ct)h(t)/(\la+\mu).$ $C=0$ is the case
%of a normal line load.

In what follows, we compute the functions needed in equations \eqref{notation} and \eqref{final}.
Equations \eqref{gtilde} and \eqref{BigG} imply
\begin{subequations}
\begin{align}
&\q{g}_1(k,l,t)=0,~\q{g}_2(k,l,t)=0,~\q{g}_3(k,l,t)=\f{\sigma_0}{\la+\mu}h(t),\\
&g^{(j)\pm}(k,l,m,t)=0,~G^{(j)}(k,l,m,t)=0,\\
&f^{(j)\pm}(k,l,m,t)=0,~F^{(j)}(k,l,m,t)=0,\\
&e^{(j)\pm}(k,l,m,t)=\f{\sigma_0}{\la+\mu} \f{1}{\pm i\om_j}(e^{\pm i\om_j t}-1),\\
&E^{(j)}(k,l,m,t)=-\f{i\sigma_0}{\la+\mu}\Big(\f{1}{\om_j^2}
-\f{\cos{(\om_j t)}}{\om_j^2}\Big),~~~j=1,2. \end{align} \end{subequations}
Thus, the known functions $N_P$, $N_Q$, $N_R$, defined in \eqref{NPNQ}, are given by
\begin{subequations}
\begin{align}
&N_P(k,l,m,t)=-m\sigma_0
\f{\la+2\mu}{\la+\mu}\Big(\f{1}{\om_1^2}-\f{\cos{(\om_1
t)}}{\om_1^2}\Big),\\
&N_Q(k,l,m,t)=l\sigma_0\f{\la+2\mu}{\la+\mu}\Big(\f{1}{\om_2^2}-\f{\cos{(\om_2
t)}}{\om_2^2}\Big),\\
&N_R(k,l,m,t)=-k\sigma_0\f{\la+2\mu}{\la+\mu}\Big(\f{1}{\om_2^2}-\f{\cos{(\om_2
t)}}{\om_2^2}\Big).
\end{align}
\end{subequations}
The above equations imply
\begin{subequations}
\begin{align}
&N_P(k,l,m,t)+N_P(k,l,-m,t)=0,\\
&N_Q(k,l,m,t)+N_Q(k,l,-m,t)=2l\sigma_0\f{\la+2\mu}{\la+\mu}\Big(\f{1}{\om_2^2}-\f{\cos{(\om_2 t)}}{\om_2^2}\Big),\\
&N_R(k,l,m,t)+N_R(k,l,-m,t)=-2k\sigma_0\f{\la+2\mu}{\la+\mu}\Big(\f{1}{\om_2^2}-\f{\cos{(\om_2 t)}}{\om_2^2}\Big),\\
&N_P(k,l,m,t)-N_P(k,l,-m,t)=-2m\sigma_0
\f{\la+2\mu}{\la+\mu}\Big(\f{1}{\om_1^2}-\f{\cos{(\om_1 t)}}{\om_1^2}\Big),\\
&N_Q(k,l,m,t)-N_Q(k,l,-m,t)=0,~~N_R(k,l,m,t)-N_R(k,l,-m,t)=0,\\
&N_P(k,l,m_{12},t)=-\f{m_{12}}{l} N_Q(k,l,m,t)=\f{m_{12}}{k} N_R(k,l,m,t),\\
&N_Q(k,l,m_{21},t)=-\f{l}{m}N_P(k,l,t),\\
&N_R(k,l,m_{21},t)=\f{k}{m}N_P(k,l,t)
\end{align}
\end{subequations}
and
\begin{subequations}\label{notation computation}
\begin{align}
&\left(
\begin{array}{llcl}
 U_2(k,l,m,t)\\
 V_2(k,l,m,t)\\
 W_2(k,l,m,t)\\
\end{array}\right)=-\f{\sigma_0}{\om_2^2}\f{\la+2\mu}{\la+\mu}C^{-1}
\left(
\begin{array}{llcl}
-m_{12}+m_{12}\cos{(\om_2 t)}\\
l-l\cos{(\om_2 t)}\\
-k+k\cos{(\om_2t)}\\
\end{array}\right),\\
&\left(
\begin{array}{llcl}
 U_1(k,l,m,t)\\
 V_1(k,l,m,t)\\
 W_1(k,l,m,t)\\
\end{array}\right)=-\f{\sigma_0}{\om_1^2}
\f{\la+2\mu}{\la+\mu}D^{-1}
\left(
\begin{array}{llcl}
m-m\cos{(\om_1 t)}\\
l-l\cos{(\om_1 t)}\\
-k+k\cos{(\om_1 t)}\\
\end{array}\right),
\end{align}
\end{subequations}
where the matrices $C$ and $D$ are defined in \eqref{Cs}. It should be noted that
the definitions of $C$ and $D$ do not depend on the particular initial-boundary values, i.e., they are
fundamental matrix-valued functions for elastodynamics in half-space. Another fundamental function appears in the alternative approach described in Appendix, which is related to Rayleigh's function, see equation (6.3).
\begin{proposition}
Let $(u,v,w)$ satisfy the Lam\'e-Navier equations \eqref{3D} in the half space
$$-\infty<x<\infty,~-\infty<y<\infty,~0<z<\infty,~t>0,$$
with the homogeneous initial conditions and Lamb's boundary conditions \eqref{Lamb_normal_point}.
A solution of this initial-boundary value problem, which decays for large $(x,y,z)$, is given by
\begin{subequations}\label{final_example}
 \begin{align}
\label{final_u_eg}
\begin{split}
&u(x,y,z,t)=\\
&\f{1}{8\pi^3}\int_{-\infty}^{\infty}dk \int_{-\infty}^{\infty} dl \int_{\gamma_{k,l}} dm~\f{e^{ikx+ily+imz}}{k^2+l^2+m^2}\\
&\Big\{ 2k[\la (k^2+l^2)+m^2(\la+2\mu)]W_1-4\mu k m^2 W_2\Big\},\\
\end{split}\\
\label{final_v_eg}
\begin{split}
& v(x,y,z,t)=\\
&\f{1}{8\pi^3}\int_{-\infty}^{\infty}dk \int_{-\infty}^{\infty} dl \int_{\gamma_{k,l}} dm~\f{e^{ikx+ily+imz}}{k^2+l^2+m^2}\\
& \Big\{ 2l [\la (k^2+l^2)+m^2(\la+2\mu)]W_1-4\mu l m^2 W_2\Big\},\\
\end{split}\\
\label{final_w_eg}
\begin{split}
& w(x,y,z,t)=\\
&\f{1}{8\pi^3}\int_{-\infty}^{\infty}dk \int_{-\infty}^{\infty} dl \int_{\gamma_{k,l}} dm~\f{e^{ikx+ily+imz}}{k^2+l^2+m^2}\\
& \Big\{4\mu k m^2 U_1+ 4\mu l m^2 V_1+2\mu k (l^2+k^2-m^2)U_2+ 2\mu l (k^2+l^2-m^2) V_2\Big\}\\
& -\f{\sigma_0}{4\pi^3}\f{\la+2\mu}{\la+\mu}\int_{-\infty}^{\infty}\!\!dk \int_{-\infty}^{\infty}\!\! dl \int_{-\infty}^{\infty}\!\! dm \f{e^{ikx+ily+imz}}{k^2+l^2+m^2}\\
& \Big\{m^2\Big[\f{1}{\om_1^2}-\f{\cos{(\om_1
t)}}{\om_1^2}\Big]+l^2\Big[\f{1}{\om_2^2}-\f{\cos{(\om_2
t)}}{\om_2^2}\Big]+k^2\Big[\f{1}{\om_2^2}-\f{\cos{(\om_2
t)}}{\om_2^2}\Big]\Big\},
\end{split}
\end{align}
\end{subequations}
where $(\om_1,\om_2)$ are defined by \eqref{dispersion_choice}, and
the known functions $(U_j,V_j,W_j)$, $j=1,2$ are computed in
\eqref{notation computation}.
\end{proposition}
Similar expressions are valid for the other Lamb's problems.

\section{Conclusions}
The main result of this paper is the derivation of equations
\eqref{final}. These equations express the displacements
$(u(x,y,z,t),v(x,y,z,t),w(x,y,z,t))$ in terms of integrals along the
real line and integrals along contours $\ga_{k,l}$ of the complex
$m-$plane; these integrals involve transforms of the given initial
and boundary data. Indeed, equations \eqref{final} involve the
functions $N_P(k,l,m,t)$, $N_Q(k,l,m,t)$ and $N_R(k,l,m,t)$, which
are defined in equations \eqref{NPNQ} in terms of the Fourier
transforms $P_j(k,l,m,t)$, $Q_j(k,l,m,t),$ $R_j(k,l,m,t)$, $j=1,2,$
of the initial data (see equations \eqref{initial}), as well as in
terms of certain known transforms $G^{(j)}(k,l,m,t),$
$F^{(j)}(k,l,m,t),$ $E^{(j)}(k,l,m,t)$, $j=1,2,$ of the boundary
data (see equations \eqref{BigG}).

The starting point of the derivations of equations \eqref{final}
is the derivation of the global relations \eqref{global_relations}. These equations are the direct
consequence of the application of the three-dimensional Fourier transform
and of the substitution in the resulting equations of the given initial
and boundary conditions. Equations \eqref{global_relations} involve the \textit{known}
functions $N_P$, $N_Q$, $N_R$, as well as the \textit{unknown} functions
$U^{(j)},V^{(j)},W^{(j)},$ $j=1,2,$ (these functions involve certain
transforms of the unknown boundary values $u(x,y,0,t)$, $v(x,y,0,t)$, $w(x,y,0,t)$ see
equations \eqref{xF} and \eqref{BigU}). The elimination of the above unknown
functions is achieved by the following steps:
1. By employing the transformation $m\rightarrow -m$, which leave the unknown
functions $(U^{(j)},V^{(j)},W^{(j)}),$ $j=1,2,$ invariant, and by utilising the analyticity
properties of these functions, we obtain the integral representations
\eqref{before_final}. These representations involve an integral along the real $k,l-axes$
and an integral along the contour $\gamma_{k,l}$ of the complex $m-$plane.
2. Using the transformations $m\rightarrow m_{21}$ and $m\rightarrow m_{12}$
which map $U^{(2)}$ to $U^{(1)}$, $V^{(2)}$ to $V^{(1)}$, $W^{(2)}$ to $W^{(1)}$ and
$U^{(1)}$ to $U^{(2)}$, $V^{(1)}$ to $V^{(2)}$, $V^{(1)}$ to $W^{(2)}$ respectively, we
express the unknown functions $(U^{(j)},V^{(j)},W^{(j)}),$ $j=1,2,$ in
terms of known functions as well as in terms of certain unknown
functions which however are analytic in a certain domain of the complex $m-$plane,
see equations \eqref{X} and \eqref{Y}. 3. Using \eqref{X} and \eqref{Y} in equations \eqref{before_final} and
employing Jordan's lemma we obtain equations \eqref{final}.

The main advantages of the new approach are the following:
\begin{enumerate}
\item{The new method provides an analytic solution of the three dimensional Lamb's problem with
arbitrary initial and boundary conditions.}

\item{This solution is expressed in terms of the given initial and
boundary data. The relevant representation is novel even for the particular case of
homogeneous initial conditions (this case has been analysed by several
authors). An alternative approach using the Laplace
transform is briefly discussed in the Appendix. By comparing
equations \eqref{final} and equation \eqref{int2}, the advantage of the new formulae
becomes clear. Actually, taking into consideration that the initial-boundary value
problems of the equations of elastodynamics are well posed for any finite $t$, the use
of the Laplace transform, which requires $t\rightarrow\infty$, is clearly inappropriate.}

\item{The new method can be employed for the solution of several related initial-boundary value problems,
including Lamb's problem of the orthotropic half space.}
\end{enumerate}

\section*{Appendix}
It is also possible to analyse Lamb's problem by using only the
transformations \eqref{superb_1} instead of using the
transformations \eqref{superb_1} \textit{and} the transformations
\eqref{superb_2}. However, in this case one \textit{cannot}
eliminate directly all unknown boundary values. Instead, one can
derive a complicated expression for the unknown boundary values in
terms of the given initial and boundary data. (This approach
is similar with the one used in \cite{Ashton_ASF} for solving
Crighton's problem).
\begin{figure}[ht]
\begin{center}
\resizebox{5.5cm}{!}{\input{Di.pstex_t}}
\end{center}
\caption{$\gamma_{k,l,2}$, the path of integration.} \label{fig:path2}
\end{figure}
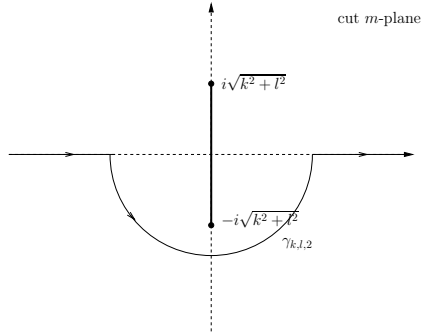
Indeed, let the contour $\gamma_{k,l,2}$ be a simple curve in the lower
half $m-$plane, determined by the requirement that it does not cross
the branch cut associated with $\om_1$ and $\om_2$, see
Fig.\ref{fig:path2}. Let $K$ be the integral operator defined by
$$(K[f])(k,l,t)=\f{1}{2\pi}\int_{\gamma_{k,l,2}}\f{f(k,l,m,t)~dm
}{m(k^2+l^2+m^2)^{1/2}},$$
for any function $f(k,l,m,t)$ with appropriate smoothness and decay. Integrating the global relations
\eqref{global_relations} along $\gamma_{k,l,2}$ we find that the functions
$$\q{h}(k,l,t)=\left(
\begin{array}{llcl}
\q{u}(k,l,t)\\
\q{v}(k,l,t)\\
\q{w}(k,l,t)\\
\end{array}\right),~~~
\q{g}(k,l,t)=\left(
\begin{array}{llcl}
\q{g}_1(k,l,t)\\
\q{g}_2(k,l,t)\\
\q{g}_3(k,l,t)\\
\end{array}\right),
$$
satisfy a system of
Volterra integral equations of the second kind:
\begin{equation}\label{int2}
\q{h}(k,l,t)=(K[N])(k,l,t)\circledast\q{g}(k,l,t)+(K[M])(k,l,t)\circledast\q{h}(k,l,t)+(K[H])(k,l,t),
\end{equation}
\begin{flushright}
$0\leq t<T, T>0;$ $k,l\in\mathbb{R},$
\end{flushright}
where $\circledast$ denotes the convolution operation with respect to
$t$,
$$N(k,l,m,t)=\left(
\begin{array}{llcl}
im\sqrt{\mu}e^{i\om_1t}                & 0 & -ik\sqrt{\mu}\f{\la+2\mu}{\mu} e^{i\om_1t}\\
0 & im\sqrt{\mu}e^{i\om_1t}  & -il\sqrt{\mu}\f{\la+2\mu}{\mu} e^{i\om_1t}\\
ik\f{\mu}{\sqrt{\la+2\mu}}e^{i\om_0t} & il\f{\mu}{\sqrt{\la+2\mu}}e^{i\om_0t}  & im\sqrt{\la+2\mu} e^{i\om_0t}\\
\end{array}\right),$$
$$M(k,l,m,t)=\left(
\begin{array}{llcl}
k^2 \sqrt{\mu}e^{i\om_1t}                & kl\sqrt{\mu} e^{i\om_1t}   & 2km\sqrt{\mu} e^{i\om_1t}\\
kl \sqrt{\mu} e^{i\om_1t}                & l^2\sqrt{\mu} e^{i\om_1t}  & 2lm\sqrt{\mu} e^{i\om_1t}\\
km\f{-2\mu}{\sqrt{\la+2\mu}}e^{i\om_1t}  & lm\f{-2\mu}{\sqrt{\la+2\mu}}e^{i\om_1t} & (k^2+l^2)\f{-\la}{\sqrt{\la+2\mu}} e^{i\om_1t}\\
\end{array}\right),$$
and the known function $H(k,l,m,t)$ is defined by
$$H(k,l,t)=\Big(i\f{i\om_1 P_0+P_1}{\la+2\mu}e^{i\om_1 t},i\f{i\om_2 Q_0+Q_1}{\mu}e^{i\om_2t},
i\f{i\om_2R_0+R_1}{\mu}e^{i\om_2t}\Big).$$
The solution of the integral equations \eqref{int2} yields the unknown transforms appearing in
\eqref{global_relations}.

Equations \eqref{int2} can be solved in closed form by using the
Laplace transform in $t$. Let ``~$\circ$~" denote the Laplace
transform with respect to $t$ and let
\begin{equation}
w_1=(p^2+(\la+2\mu)(k^2+l^2))^{\f{1}{2}},~w_2=(p^2+\mu
(k^2+l^2))^{\f{1}{2}},
\end{equation}
\begin{equation}
\Delta=\f{(p^2+2\mu (k^2+l^2))^2}{w_1^2w_2^2}-4\sqrt{\f{\mu}{\la+2\mu}}\f{\mu (k^2+l^2)}{w_1 w_2}.
\end{equation}
The solution of \eqref{int2} is given by
\begin{equation}\label{Spectral_NDC}
\mathring{\q{h}}=\f{1}{\Delta}(I-\mathring{K[M]})^{*T}(\mathring{K[N]}\mathring{\q{g}}+\mathring{K[H]}),
\end{equation}
where ``$~\ast~$" denotes the adjoint operation, and the functions
$\mathring{K[M]}$, $\mathring{K[N]}$ which are independent of the
initial and boundary conditions, are given by
\begin{equation}
 \mathring{K[M]}=\left(\begin{array}{llcl}
\f{-k^2\mu}{w_2^2}                & \f{-kl \mu}{w_2^2}  & \f{2ik\sqrt{\mu}}{w_2}\\
\f{-kl\mu}{w_2^2}                & \f{-l^2 \mu}{w_2^2}  & \f{2il\sqrt{\mu}}{w_1}\\
\sqrt{\f{\mu}{\la+2\mu}}\f{-2ik\sqrt{\mu}}{w_1}& \sqrt{\f{\mu}{\la+2\mu}}\f{-2il\sqrt{\mu}}{w_1} & \f{-\la (k^2+l^2)}{w_1}\\
\end{array}\right)
\end{equation}
and
\begin{equation}
\mathring{K[N]}=\left(\begin{array}{llcl}
-\f{\sqrt{\mu}}{w_2} & 0               & i\f{(\la+2\mu)k}{w_2^2}\\
0   &   -\f{\sqrt{\mu}}{w_2}               & i\f{(\la+2\mu)k}{w_2^2}\\
\f{-ik\mu}{w_1^2} &\f{-il\mu}{w_1^2} & \f{-\sqrt{\la+2\mu}}{w_1}\\
\end{array}\right).
\end{equation}
The zeros of $\Delta$ coincide with the zeros of Rayleigh's
function. When $\mu/\la>0.906$, the known
transforms $\mathring{\q{g}}(k,l,p)$, for each fixed $k,l$, do \emph{not} have poles with
positive real parts.

In the particular case of the problem of homogeneous initial
condition and the normal point load boundary condition, equation
\eqref{Spectral_NDC} reduces to the classic Lamb's solution
\cite{HL}\cite{JM}.

\paragraph{Acknowledgements}
We would like to thank Ishan Sharma for drawing our attentions on
Lamb's problems. We are grateful to Michail Dimakos, Dionyssis
Mantzavinos and Anthony Ashton for several valuable observations. Di
Yang would like to thank Professor Youjin Zhang of Tsinghua
University for his advise and for several discussions, as well as
the China Scholarship Council for supporting him for a joint PhD
study at the University of Cambridge. A. S. Fokas would like to
thank Y. Antipov for his suggestions and also acknowledges the
generous support of the Guggenheim Foundation, USA. The work of Yang
is partially supported by the National Basic Research Program of
China (973 Program) No.2007CB814800. 

\bibliographystyle{amsplain}

\end{document}

%% file: l12.pstex_t
\begin{picture}(0,0)%
\includegraphics{l12.pstex}%
\end{picture}%
\setlength{\unitlength}{3947sp}%
\begingroup\makeatletter\ifx\SetFigFont\undefined%
\gdef\SetFigFont#1#2#3#4#5{%
  \reset@font\fontsize{#1}{#2pt}%
  \fontfamily{#3}\fontseries{#4}\fontshape{#5}%
  \selectfont}%
\fi\endgroup%
\begin{picture}(6024,4149)(2614,-5923)
\put(3376,-3661){\makebox(0,0)[lb]{\smash{{\SetFigFont{14}{16.8}{\familydefault}{\mddefault}{\updefault}{\color[rgb]{0,0,0}$-\sqrt{k^2+l^2}\sqrt{\f{\la+\mu}{\mu}}$}%
}}}}
\put(7051,-2161){\makebox(0,0)[lb]{\smash{{\SetFigFont{14}{16.8}{\familydefault}{\mddefault}{\updefault}{\color[rgb]{0,0,0}cut $m$-plane}%
}}}}
\put(5776,-3661){\makebox(0,0)[lb]{\smash{{\SetFigFont{14}{16.8}{\familydefault}{\mddefault}{\updefault}{\color[rgb]{0,0,0}$\sqrt{k^2+l^2}\sqrt{\f{\la+\mu}{\mu}}$}%
}}}}
\end{picture}%

%% file: l21.pstex_t
\begin{picture}(0,0)%
\includegraphics{l21.pstex}%
\end{picture}%
\setlength{\unitlength}{3947sp}%
\begingroup\makeatletter\ifx\SetFigFont\undefined%
\gdef\SetFigFont#1#2#3#4#5{%
  \reset@font\fontsize{#1}{#2pt}%
  \fontfamily{#3}\fontseries{#4}\fontshape{#5}%
  \selectfont}%
\fi\endgroup%
\begin{picture}(6024,4149)(2614,-5923)
\put(5026,-5311){\makebox(0,0)[lb]{\smash{{\SetFigFont{14}{16.8}{\familydefault}{\mddefault}{\updefault}{\color[rgb]{0,0,0}$-i\sqrt{k^2+l^2}\sqrt{\f{\la+\mu}{\la+2\mu}}$}%
}}}}
\put(5026,-2686){\makebox(0,0)[lb]{\smash{{\SetFigFont{14}{16.8}{\familydefault}{\mddefault}{\updefault}{\color[rgb]{0,0,0}$i\sqrt{k^2+l^2}\sqrt{\f{\la+\mu}{\la+2\mu}}$}%
}}}}
\put(7051,-2161){\makebox(0,0)[lb]{\smash{{\SetFigFont{14}{16.8}{\familydefault}{\mddefault}{\updefault}{\color[rgb]{0,0,0}cut $m$-plane}%
}}}}
\end{picture}%

%% file: Newversion_1.pstex_t
\begin{picture}(0,0)%
\includegraphics{Newversion_1.pstex}%
\end{picture}%
\setlength{\unitlength}{3947sp}%
\begingroup\makeatletter\ifx\SetFigFont\undefined%
\gdef\SetFigFont#1#2#3#4#5{%
  \reset@font\fontsize{#1}{#2pt}%
  \fontfamily{#3}\fontseries{#4}\fontshape{#5}%
  \selectfont}%
\fi\endgroup%
\begin{picture}(6058,4149)(2239,-5923)
\put(7051,-2161){\makebox(0,0)[lb]{\smash{{\SetFigFont{14}{16.8}{\familydefault}{\mddefault}{\updefault}{\color[rgb]{0,0,0}cut $m$-plane}%
}}}}
\put(3451,-2761){\makebox(0,0)[lb]{\smash{{\SetFigFont{14}{16.8}{\familydefault}{\mddefault}{\updefault}{\color[rgb]{0,0,0}$\gamma_{k,l}$}%
}}}}
\end{picture}%

%% file: Di.pstex_t
\begin{picture}(0,0)%
\includegraphics{Di.pstex}%
\end{picture}%
\setlength{\unitlength}{3947sp}%
\begingroup\makeatletter\ifx\SetFigFont\undefined%
\gdef\SetFigFont#1#2#3#4#5{%
  \reset@font\fontsize{#1}{#2pt}%
  \fontfamily{#3}\fontseries{#4}\fontshape{#5}%
  \selectfont}%
\fi\endgroup%
\begin{picture}(6133,4899)(2389,-6898)
\put(7276,-2311){\makebox(0,0)[lb]{\smash{{\SetFigFont{14}{16.8}{\familydefault}{\mddefault}{\updefault}{\color[rgb]{0,0,0}cut $m$-plane}%
}}}}
\put(6451,-5611){\makebox(0,0)[lb]{\smash{{\SetFigFont{14}{16.8}{\familydefault}{\mddefault}{\updefault}{\color[rgb]{0,0,0}$\gamma_{k,l,2}$}%
}}}}
\put(5551,-3286){\makebox(0,0)[lb]{\smash{{\SetFigFont{14}{16.8}{\familydefault}{\mddefault}{\updefault}{\color[rgb]{0,0,0}$i\sqrt{k^2+l^2}$}%
}}}}
\put(5551,-5311){\makebox(0,0)[lb]{\smash{{\SetFigFont{14}{16.8}{\familydefault}{\mddefault}{\updefault}{\color[rgb]{0,0,0}$-i\sqrt{k^2+l^2}$}%
}}}}
\end{picture}%